\documentclass[12pt,reqno]{amsart}                  
\usepackage{amssymb}
\usepackage{mathrsfs}

\newtheorem{theorem}{Theorem}[section]
\newtheorem{lemma}[theorem]{Lemma}

\theoremstyle{remark}
\newtheorem{remark}[theorem]{Remark}
\begin{document}
\title[Lower bounds of martingale measure densities]{Lower bounds of martingale measure densities in the Dalang-Morton-Willinger theorem}

\author{Dmitry B. Rokhlin}

\address{D.B. Rokhlin,
Faculty of Mathematics, Mechanics and Computer Sciences,
              Southern Federal University, 
Mil'chakova str., 8a, 344090, Rostov-on-Don, Russia}          
\email{rokhlin@math.rsu.ru}                   
\newcommand{\ri}{{\rm ri\,}}
\newcommand{\gr}{{\rm gr\,}}
\newcommand{\cl}{{\rm cl}}
\newcommand{\epi}{{\rm epi\,}}
\newcommand{\cone}{{\rm cone\,}}
\newcommand{\intern}{{\rm int\,}}
\newcommand{\conv}{{\rm conv\,}}
\newcommand{\dom}{{\rm dom\,}}
\newcommand{\supp}{{\rm supp\,}}
\newcommand{\esssup}{{\rm ess\,sup}}
\renewcommand{\theequation}{\thesection.\arabic{equation}}
\thispagestyle{empty}

\begin{abstract} For a $d$-dimensional stochastic process $(S_n)_{n=0}^N$ we obtain criteria for the existence of an equivalent martingale measure, whose density $z$, up to a normalizing constant, is bounded from below by a given random variable $f$. We consider the case of one-period model ($N=1$) under the assumptions $S\in L^p$; $f,z\in L^q$, $1/p+1/q=1$, where $p\in [1,\infty]$, and the case of $N$-period model for $p=\infty$. The mentioned criteria are expressed in terms of the conditional distributions of the increments of $S$, as well as in terms of the boundedness from above of an utility function related to some optimal investment problem under the loss constraints. Several examples are presented.   
\end{abstract}

\maketitle

\section*{Introduction} \label{intro}
\setcounter{equation}{0}
Let $(\Omega,\mathscr F,\mathsf P)$ be a probability space, endowed with a discrete-time filtration $\mathbb F=(\mathscr F_n)_{n=0}^N$, $\mathscr F_N=\mathscr F$. Consider a $d$-dimensional stochastic process $S=(S_n)_{n=0}^N$, adapted to the filtration $\mathbb F$, and a $d$-dimensional $\mathbb F$-predictable process $\gamma=(\gamma_n)_{n=1}^N$. In the customary securities market model $S^i_n$ describes the discounted price of $i$th stock and $\gamma^i_n$ corresponds to the number of stock units in investor's portfolio at time moment $n$. The gain process is given by
\begin{equation}
G_n^\gamma=\sum_{k=1}^n(\gamma_k,\Delta S_k),\ \ \ \Delta S_k=S_k-S_{k-1},\ \ n=1,\dots,N,
\end{equation}
where $(a,b)$ is the scalar product of $a,b\in\mathbb R^d$. 

Let's recall the classical Dalang-Morton-Willinger theorem \cite{DalMorWil90}, \cite{Shi98} (ch.V, \S 2e).
As usual, we say that the \emph{No Arbitrage (NA)} condition is satisfied if the inequality $G_N^\gamma\ge 0$ a.s. (with respect to the measure $\mathsf P$) implies that $G_N^\gamma=0$ a.s. A probability measure $\mathsf Q$ on $\mathscr F$ is called a \emph{martingale measure} if the process $S$ is a $\mathsf Q$-martingale. The measures $\mathsf P$ and $\mathsf Q$ are called \emph{equivalent} if their null sets are the same. Denote by $\varkappa_{n-1}(\omega)$ the \emph{support of the regular conditional distribution} $\mathsf P_{n-1}(\omega,dx)$ of the random vector $\Delta S_n$ with respect to $\mathscr F_{n-1}$.
\begin{theorem}[Dalang-Morton-Willinger]
The following conditions are equivalent:
\begin{itemize}
\item[(i)] NA;
\item[(ii)] there exists an equivalent to $\mathsf P$ martingale measure $\mathsf Q$ with a.s. bounded density $z=d\mathsf Q/d\mathsf P$;
\item[(iii)] the relative interior of the convex hull of $\varkappa_{n-1}$ contains the origin a.s., $n=1,\dots,N$.
\end{itemize}
\end{theorem}

The question concerning the existence of an equivalent martingale measure $\mathsf Q$, whose density $z$ satisfies the the lower bound $z\ge c$ (where $c$ is a positive constant) was posed in \cite{RasSte05} (Remark 7.5), \cite{DelSch06} (Remark 6.5.2). In general, the answer to this question is negative. An evident necessary condition is the integrability of $S$ with respect to $\mathsf P$. Moreover, the example of \cite{DelSch06} shows that a measure $\mathsf Q$ with the above properties need not exist even for a uniformly bounded process $S$. A sufficient condition was obtained in \cite{RasSte05}. In particular it is satisfied for a process $S$ with independent increments, if the random vectors $\Delta S_n$ have finite moments of all orders. 

Following \cite{RokSch06}, let us formulate the problem concerning the existence of an equivalent martingale measure, whose density (up to a normalization constant) is bounded from below by a random variable $f$, in a more general context. Denote by $\mathsf E X$ the expectation with respect to $\mathsf P$, by $L^p=L^p(\mathscr F)=L^p(\Omega,\mathscr F,\mathsf P)$, $p\in [1,\infty)$ the Banach spaces of equivalence classes of $\mathscr F$-measurable functions with the norms $\|X\|_p=\mathsf E |X|^p$ and by $L^\infty$ the Banach space of essentially bounded functions with the norm $\|X\|_\infty=\esssup|X|$. The cone $L^p_+$ of non-negative elements induces the partial order on $L^p$ .

Consider the subspace $K\subset L^p$, $p\in [1,\infty)$ of investor's gains (discounted wealth increments). Denote by $q$ the conjugate exponent, that is, $1/p+1/q=1$. The condition $K\cap L^p_+=\{0\}$ corresponds to NA. An element $f\in L^q_+$ induces the functional on $L^p$ by the formula $\langle X,f\rangle=\mathsf E(Xf)$, $X\in L^p$. It turns out that the  existence of an element $g$, satisfying the conditions
\begin{equation}
\langle X,g\rangle=0,\ \  X\in K;\ \ \ g\ge f,\ \ g\in L^q
\end{equation}
is equivalent to the boundedness of $f$ form above on a certain subset $K_1$ of the subspace $K$:
\begin{equation}
v_p:=\sup_{X\in K_1}\langle X,f\rangle<\infty,\ \ \ K_1=\{X\in K:\|X^-\|_p\le 1\},
\end{equation}
where $X^-=\max\{-X,0\}$. For $p=\infty$, $q=1$ this statement is not true in general, see \cite{RokSch06}, Examples 1 and 3. It becomes true under the assumption that $f$ is bounded from above on the subset $\{X\in K: X^-\in V\},$
where $V$ is a neighborhood of zero in the Mackey topology $\tau(L^\infty,L^1)$, or if $L^1$ is replaced by the topological dual space $(L^\infty)^*$ of $L^\infty$. These results are contained in Theorem 1 of \cite{RokSch06}.

It should be mentioned that the problems, equivalent to (0.3) when $f=1$, were considered in the recent paper\cite{Lei08}. From the financial point of view they correspond to the maximization of expected gain under the loss constraint, if the loss value is measured either by $p$th moment $\mathsf E |X^-|^p$ for $p\in [1,\infty)$ or by $\esssup|X^-|$ for $p=\infty$. The equivalence of (0.2) and (0.3) for $p\in (1,\infty)$ follows from the results of the cited paper as well (\cite{Lei08}, Theorem 4.1). Unfortunately, the related statement for $p=\infty$ (\cite{Lei08}, Theorem 6.1) is incorrect: a counterexample is, in fact, contained in \cite{RokSch06} (Example 3) and its another version is given below (Example 5.4).  

Turning back to the finite securities market model, assume that $S\in L^p$ and denote by $K$ the set of random variables $G_N^\gamma$, where $\gamma$ is a \emph{bounded} predictable process. Then the elements $g$, satisfying (0.2), up to a normalization constant, coincide with the $\mathsf P$-densities of martingale measures: $d\mathsf Q/d\mathsf P=g/\mathsf E g$. 

The aim of the present paper is to establish effective criteria for the fulfilment of (0.2), (0.3) for a market model with \emph{finite discrete time} and a \emph{finite collection of stocks}. Such criteria, expressed in terms of the regular conditional distributions of the increments $\Delta S_n$, are obtained for a one-period model under the assumptions $S\in L^p$, $f,g\in L^q$, $p\in [1,\infty]$ (Theorem 1.3), as well as for $N$-period model in the case $p=\infty$ (Theorem 4.1). These results show also that in the case under consideration the equivalence of (0.2) and (0.3) for $p=\infty$ is nevertheless true! Thereby, we give the negative answer to the question, raised in the end of the paper \cite{RokSch06}. 

In the last part of the paper we give some examples, illustrating the effectiveness of the obtained criteria, and a counterexample to the mentioned statement of \cite{Lei08}. Also, it is interesting to note that the case $p=1$ of Theorem 1.2 leads to a new proof of the key implication (iii) $\Longrightarrow$ (ii) of the Dalang-Morton-Willinger theorem (Remark 1.5). 

\section{One-period model}
\label{sec:1}
\setcounter{equation}{0}
Let $(\Omega,\mathscr F,\mathsf P)$ be a probability space and let $\mathscr H$ be a sub-$\sigma$-algebra of $\mathscr F$. A set-valued mapping $F$, assigning some set $F(\omega)\subset\mathbb R^d$ to each $\omega\in\Omega$, is called  $\mathscr H$-\emph{measurable}, if $\{\omega: F(\omega)\cap V\neq\varnothing\}\in\mathscr H$ for any open set $V\subset\mathbb R^d$. A function $\eta:\Omega\mapsto\mathbb R^d$ is called a \emph{selector} of $F$, if $\eta(\omega)\in F(\omega)$ for all $\omega\in\dom F:=\{\omega':F(\omega')\neq\varnothing\}$.
An $\mathscr H$-measurable set-valued mapping $F$ with non-empty closed values $F(\omega)$ is measurable if and only if there exists a sequence $(\eta_i)_{i=1}^\infty$ of $\mathscr H$-measurable selectors of $F$ such that the sets $\{\eta_i(\omega)\}_{i=1}^\infty$ are dense in $F(\omega)$ for all $\omega$ (\cite{Rockaf76}, Theorem 1B). Such a sequence is called a \emph{Castaing representation} of $F$. 

Denote by $\mathscr B(\mathbb R^d)$ the Borel $\sigma$-algebra of $\mathbb R^d$. A function $\varphi:\Omega\times\mathbb R^d\mapsto\mathbb R$ is called a \emph{Carath\'eodory function} if (a) $\varphi(\cdot,x):\Omega\mapsto\mathbb R$ is $(\mathscr H,\mathscr B(\mathbb R))$-measurable for all $x\in\mathbb R^d$, (b) $\varphi(\omega,\cdot):\mathbb R^d\mapsto\mathbb R$ is continuous for all $\omega\in\Omega$. 

Denote by $L^p(\mathscr H,F)$, $1\le p<\infty$ the set of equivalence classes of $\mathscr H$-measurable functions $\eta$ satisfying the conditions $\int |\eta|^p\,d\mathsf P<\infty$, $\eta\in F$ a.s., where $|x|=(x,x)^{1/2}$. We introduce also the sets of equivalence classes  of essentially bounded functions $L^\infty(\mathscr H,F)$ and of all $\mathscr H$-measurable functions $L^0(\mathscr H,F)$, taking values in $F$ a.s. In accordance with the above notation we put $L^p(\mathscr H)=L^p(\mathscr H,\mathbb R)$. By $L^p_+(\mathscr H)$ and $L^p_{++}(\mathscr H)$ we denote the sets of non-negative and strictly positive elements of $L^p(\mathscr H)$ respectively. Let $\|X\|_p$ be the norm of an element $X$ of the Banach space $L^p(\mathscr H)$, $1\le p\le\infty$. 

The completion of the $\sigma$-algebra $\mathscr H$ with respect to the measure $\mathsf P$ is denoted by $\mathscr H^\mathsf P$. Note that $L^p(\mathscr H^\mathsf P)=L^p(\mathscr H)$ in the sense that any $\mathscr H^\mathsf P$-measurable function possesses an $\mathscr H$-measurable modification.

In the sequel we use the customary notation of convex analysis for the polar $A^\circ=\{x\in\mathbb R^d: (x,y)\le 1, y\in A\}$ of a set $A\subset\mathbb R^d$ and also for its Minkowski function and the support function:
$$ \mu(x|A)=\inf\{\lambda>0:x\in\lambda A\},\ \ \ s(x|A)=\sup_{y\in A}(x,y).$$ 
Denote by $\conv A$, $\ri A$ the convex hull and the relative interior of $A$.

Consider the one-period model (0.1) (that is, $N=1$). Put $\xi=\Delta S_1$, $\mathscr H=\mathscr F_0$. Let $\mathsf P_\xi(\omega,dx)$ be the regular conditional distribution of $\xi$ with respect to $\mathscr H$ and let $\varkappa_\xi(\omega)$ be the support of the measure $\mathsf P_\xi(\omega,\cdot)$. By $D_\xi(\omega)\subset\mathbb R^d$ we denote the linear span of $\varkappa_\xi(\omega)$.
Define the functions
\begin{eqnarray*}
\psi_p(\omega,h) &=& \left(\int_{\mathbb R^d} [(h,x)^{-}]^p\,\mathsf P_\xi(\omega,dx)\right)^{1/p},\ p\in[1,\infty);\\
\psi_\infty(\omega,h) &=&s(-h|\varkappa_\xi(\omega))
\end{eqnarray*}
from $\Omega\times\mathbb R^d$ to $[0,\infty]$, and the set-valued mappings
 \begin{equation}     
 \omega\mapsto T_p(\omega)=\{h\in D_\xi(\omega):\psi_p(\omega,h)\le 1\}.
 \end{equation}      

\begin{lemma}
Assume that $0\in\ri(\conv\varkappa_\xi(\omega))$ a.s. Then $T_p$ is an $\mathscr H^\mathsf P$-measurable set-valued mapping with a.s. compact values, $p\in [1,\infty]$. 
\end{lemma}
{\sc Proof.} The set-valued mapping $\omega\mapsto\varkappa_\xi(\omega)$ is $\mathscr H$-measurable:
$$ \{\omega:\varkappa_\xi(\omega)\cap V\neq\varnothing\}=\{\omega:\mathsf P_\xi(\omega,V)>0\}\in\mathscr H$$
for any open set $V\subset\mathbb R^d$. Its values $\varkappa_\xi(\omega)$ are closed. It follows from the formula
$$ \psi_\infty(\omega,h)=\sup_{i\ge 1}(-h,\eta_i(\omega)),$$
where $(\eta_i)_{i=1}^\infty$ is a Castaing representation of $\varkappa_\xi$, that the function $\omega\mapsto\psi_\infty(\omega,h)$ is $\mathscr H$-measurable. The same property of $\psi_p$ for $p\in[1,\infty)$ is evident.

Put $\Omega_p=\{\omega:\int |x|^p\,d\mathsf P_\xi(\omega,dx)<\infty\}$ for $p\in [1,\infty)$ and let $\Omega_\infty$ be  the set of $\omega$, for which the set $\varkappa_\xi(\omega)$ is compact. Note that $\Omega_\infty=\{\omega:\sup_{h\in\mathbb D}\psi_\infty(\omega,h)<\infty\}$, where $\mathbb D$ is a countable dense subset of $\mathbb R^d$. Consequently, $\Omega_p\in\mathscr H$, $p\in [1,\infty]$ and $\mathsf P(\Omega_p)=1$. Put $\Omega_p'=\Omega_p\cap\{\omega:0\in\ri(\conv\varkappa_\xi(\omega))\}$. Clearly, $\Omega_p'\in\mathscr H^\mathsf P$ and $\mathsf P(\Omega'_p)=1$.

Assume that $\omega\in\Omega_p'$. It follows from continuity of $\psi_p$ with respect to $h$ that the set $T_p(\omega)$ is closed. From the codition $0\in\ri(\conv\varkappa_\xi(\omega))$ we see that for $h\in D_\xi(\omega)\backslash 0$ the set $\varkappa_\xi(\omega)$ is not contained in the half-space $\{x\in D_\xi(\omega):(h,x)\ge 0\}$. Therefore, $\psi_p(\omega,h)>0$, $p\in [1,\infty]$ and the set $T_p(\omega)$ is compact, because 
$\psi_p(\omega,h)\to \infty$ when $|h|\to\infty$, $h\in D_\xi(\omega)$. 

Consider the trace of the $\sigma$-algebra $\mathscr H$ on $\Omega_p'$: $\mathscr H_p=\{A\cap\Omega_p':A\in\mathscr H\}$. To complete the proof it is sufficient to check that the set-valued mappings $\omega\mapsto T_p(\omega)$, $\omega\in\Omega_p'$ are $\mathscr H_p$-measurable. We make use of the representation
$T_p(\omega)=\{h\in\mathbb R^d:\psi_p(\omega,h)\le 1\}\cap D_\xi(\omega)$ and the fact that $\psi_p:(\Omega_p'\times\mathbb R^d,\mathscr H_p\otimes\mathscr B(\mathbb R^d))\mapsto [0,\infty)$ are Carath\'eodory functions. The measurability of the each of set-valued mappings, whose intersection is $T_p$, follows from Corollary 1Q and Proposition 1H of \cite{Rockaf76}, and the measurability of $T_p$ is implied by Theorem 1M of the same paper. $\square$ 

Let us recall the ''measurable maximum theorem'' (\cite{AliBor07}, Theorem 18.19).
\begin{lemma}  \label{lem1.2}
Let $F$ be an $\mathscr H$-measurable set-valued mapping with non-empty compact values $F(\omega)\subset\mathbb R^d$, and let $\varphi:\Omega\times \mathbb R^d\mapsto\mathbb R$ be a Carath\'eodory function. Put
$$ m(\omega)=\max_{x\in F(\omega)} \varphi(\omega,x),\ \ \ G(\omega)=\{x\in F(\omega):\varphi(\omega,x)=m(\omega)\}.$$
Then (a) the function $m$ and the set-valued mapping $G$ are $\mathscr H$-measurable; 
(b) there exists an $\mathscr H$-measurable selector $\eta^*$ of $G$.
\end{lemma}

Our first main result is the following.
\begin{theorem} \label{teo1.3}
Let $\xi\in L^p(\mathscr F,\mathbb R^d)$, $f\in L^q_+(\mathscr F)$, where $p\in [1,\infty]$ and $1/p+1/q=1$. If  $0\in\ri(\conv\varkappa_\xi)$ a.s., then the following conditions are equivalent:
\begin{itemize}
\item[(i)] $v_p:=\sup\{\mathsf E( fX): \|X^-\|_p\le 1,\ X\in K\}<\infty$, where 
    \begin{equation}
K=\{(\gamma,\xi):\gamma\in L^\infty(\mathscr H,D_\xi)\};
   \end{equation}     
\item[(ii)] there exists a random variable $g\in L^q(\mathscr F)$, satisfying the conditions
    \begin{equation}
           \mathsf E(g\xi|\mathscr H)=0,\ \ \ g\ge f;
    \end{equation}       
\item[(iii)] $s(a|T_p)\in L^q(\mathscr H)$, where $a=\mathsf E(f\xi|\mathscr H)$ and $T_p$ is defined by the formula (1.1).
\end{itemize}
\end{theorem}

Let us make some remarks before the proof of this theorem (sect. 2 and 3). 

\begin{remark} \label{rem:1.4} If $0\in\ri(\conv\varkappa_\xi)$ and $\xi\in L^1(\mathscr F,\mathbb R^d)$ does not depend on $\mathscr H$, then there exists $g\in L^\infty(\mathscr F)$:
$$ \mathsf E(g\xi|\mathscr H)=0,\ \ \ g\ge 1.$$
Actually, in this case $s(a|T_1)=s(\mathsf E\xi|T_1)$ does not depend on $\omega$ and thus belongs to $L^\infty(\mathscr H)$. 
\end{remark}

\begin{remark} \label{rem:1.5} If $0\in\ri(\conv\varkappa_\xi)$ and $\xi\in L^1(\mathscr F,\mathbb R^d)$, then there exists $g\in L^\infty_{++}(\mathscr F)$:
$$ \mathsf E(g\xi|\mathscr H)=0.$$
To prove this statement it is sufficient to note that there exists an $\mathscr H$-measurable function $f\in L^\infty_{++}(\mathscr H)$ such that 
$$ s(\mathsf E(f\xi|\mathscr H)|T_1)=s(\mathsf E(\xi|\mathscr H)|T_1) f\in L^\infty(\mathscr H).$$
A function $g\in L^\infty(\mathscr F)$, satisfying (1.3), is the desired one. 

In fact, this proves the implication (iii) $\Longrightarrow$ (ii) of Theorem 0.1 for $N=1$ and $S\in L^1$. As is known, this is the key point of the proof of the Dalang-Morton-Willinger theorem.
\end{remark}

\begin{remark} \label{rem:1.6} Note that
\begin{equation}
 a=\mathsf E(\mathsf E(f|\mathscr H\vee\sigma(\xi))\xi|\mathscr H)=\int b(\omega,x)x\,\mathsf P_\xi(\omega,dx) \in D_\xi(\omega)\ \ \textrm{a.s.}
\end{equation} 
The existence of an $\mathscr H\otimes\mathscr B(\mathbb R^d)$-measurable function $b(\omega,x)$, satisfying the condition $\mathsf E(f|\mathscr H\vee\sigma(\xi))=b(\omega,\xi)$, follows from the fact that the $\sigma$-algebra $\mathscr H\vee\sigma(\xi)$ is generated by the mapping $\omega\mapsto(\omega,\xi(\omega))$ from $\Omega$ to the measurable space $(\Omega \otimes\mathbb R^d,\mathscr H\otimes\mathscr B(\mathbb R^d))$.  
\end{remark}

\begin{remark} \label{rem:1.7} We have the following convenient representation of the random variable $s(a(\omega)|T_\infty(\omega))$ for $a\in D_\xi$ a.s.: 
$$ s(a|T_\infty)=\sup\{(h,a): -h\in D_\xi\cap \varkappa_\xi^\circ\}=s(-a|\varkappa_\xi^\circ)=\mu(-a|\conv\varkappa_\xi)\ \textrm{a.s.}$$
It the last equality we have used the formula
$$ \mu(x|A^\circ)=\inf\{\lambda>0:\lambda^{-1} x\in A^\circ\}=\inf\{\lambda>0:s(\lambda^{-1} x|A)\le 1\}=
s(x|A),$$
which is true under the assumption $0\in A$. We have also used the bipolar theorem: $A^{\circ\circ}=\cl(\conv A)$ and the compactness property of the convex hull of a compact set. 
\end{remark}

\section{Proof of Theorem 1 for $p\in [1,\infty)$}
\setcounter{equation}{0}
Denote by $U^p$ the unit ball of the space $L^p=L^p(\Omega,\mathscr F,\mathsf P)$ and put
$U_+^p=\{X\in L^p_+:X\in U^p\}$. 
\begin{lemma} \label{lem2.1}
For any element $X\in L^p$, $p\in [1,\infty]$ we have
$$\|X^+\|_p=\sup\{\langle X,z\rangle:z\in U^q_+\},\ \ \frac{1}{p}+\frac{1}{q}=1.$$
\end{lemma}
{\sc Proof.} Consider the elements 
$$\zeta_q=\frac{(X^+)^{p/q} }{\|X^+\|_p^{p/q}}\in U^q_+,\ q\in (1,\infty);\ \ \zeta_\infty=I_{\{X\ge 0\}}\in U_+^\infty; \ \ \zeta_1^n=\frac{I_{A_n}}{\mathsf P(A_n)}\in U^1_+,$$
where $A_n=\{\omega:X(\omega)\ge\|X^+\|_\infty-1/n\}$. If $X\in L^p$ and $q$ is the conjugate exponent, then 
$$\langle X,\zeta_q\rangle=\|X^+\|_p,\ q\in (1,\infty]; \ \ \langle X,\zeta_1^n\rangle\ge \|X^+\|_\infty-\frac{1}{n}.$$
On the other hand,
$$ \langle X,z\rangle\le\langle X^+,z\rangle\le\|X^+\|_p,\ \ z\in U^q_+.\ \ \square$$

Though the next result follows from Theorem 1 of \cite{RokSch06}, it seems convenient to give its direct proof. The idea of this proof is contained also in the paper \cite{Rok05} (Lemma 2.5).

Recall that the closure of a convex set $A\subset L^p$, $p\in [1,\infty)$ in the weak topology $\sigma (L^p,L^q)$, $1/p+1/q=1$ coincides with its norm closure in $L^p$.  
\begin{lemma} \label{lem2.2}
For a subspace $K\subset L^p$, $p\in [1,\infty)$ and an element $f\in L^q_+$, $1/p+1/q=1$ the following conditions are equivalent:
\begin{itemize}
\item[(a)] $\sup_{X\in K_1}\langle X,f\rangle<\infty$, where $K_1=\{X\in K:\|X^-\|_p\le 1\};$
\item[(b)] there exists $g\in L^q$, satisfying the conditions 
 \begin{equation}     
 \langle X,g\rangle=0,\ \ X\in K; \ \ g\ge f.
 \end{equation}     
\end{itemize}
\end{lemma}
{\sc Proof.} (b) $\Longrightarrow$ (a). If $X\in K_1$ then
$$ \langle X,f\rangle=\langle X,g\rangle+\langle X,f-g\rangle=-\langle X,g-f\rangle\le \langle X^-,g-f\rangle\le\|g-f\|_q.$$

(a) $\Longrightarrow$ (b). Put $\lambda=\sup_{X\in K_1}\langle X,f\rangle$. 
If the assertion (b) is false then
$$ (f+\lambda U^q_+)\cap K^\circ=\varnothing,\ \ \ K^\circ=\{z\in L^q:\langle X,z\rangle\le 0,\ X\in K\}.$$
By applying the separation theorem (\cite{AliBor07}, Theorem 5.79) to the $\sigma (L^q,L^p)$-compact set $f+\lambda U^q_+$ and to the $\sigma (L^q,L^p)$-closed set $K^\circ$, we conclude that there exists $Y\in L^p$ such that
$$ \sup_{z\in K^\circ}\langle Y,z\rangle<\inf\{\langle Y,\zeta\rangle:\zeta\in f+\lambda U^q_+\}.$$
Since $K$ is a subspace it follows that $\langle Y,z\rangle=0$, $z\in K^\circ$ and $Y\in K^{\circ\circ}=\cl_p K$ by the bipolar theorem (\cite{AliBor07}, Theorem 5.103), where $\cl_p K$ is the closure of $K$ in the norm topology of $L^p$. Moreover,
\begin{eqnarray}
\langle Y,f\rangle+\lambda\inf\{\langle Y,\eta\rangle:\eta\in U^q_+\}>0.
\end{eqnarray}

By Lemma 2.1 we have
\begin{eqnarray}
\inf\{\langle Y,\eta\rangle:\eta\in U^q_+\}=-\sup\{\langle -Y,\eta\rangle:\eta\in U^q_+\}=-\|Y^-\|_p.
\end{eqnarray}
If $Y^-=0$ then $\langle Y,f\rangle>0$ and $\alpha Y\in L^p_+\cap\cl_p K$ for any $\alpha>0$. Hence, the functional $X\mapsto \langle X,f\rangle$ is unbounded from above on the ray $\{\alpha Y:\alpha>0\}$, which lies in the set
$$\cl_p K_1\supset\cl_p\biggr(\{X:\|X^-\|_p<1\}\cap K\biggl)\supset \{X:\|X^-\|_p<1\}\cap \cl_p K.$$
Here we have used the elementary inclusion
$ \cl_p(A\cap B)\supset A\cap\cl_p B,$
which holds true when the set $A$ is open (\cite{Bur68}, chap.1, \S1, Proposition 5).

Thus, $\|Y^-\|_p>0$. It follows from (2.2), (2.3) that 
$$ \langle Y/\|Y^-\|_p,f\rangle>\lambda$$
in contradiction with the definition of $\lambda$ since $Y/\|Y^-\|_p\in K_1$. $\square$

Lemma 2.2 implies that the conditions (i) and (ii) of Theorem 1.3 are equivalent. Indeed, for the subspace (1.2) condition $\langle X,g\rangle=0$, $X\in K$ means that
\begin{equation}
\mathsf E[g(\gamma,\xi)]=\mathsf E(\gamma,\mathsf E(g\xi|\mathscr H))=0,\ \ \ \gamma\in L^\infty(\mathscr H,D_\xi).
\end{equation}
In turn, (2.4) is reduced to the equality $\mathsf E(g\xi|\mathscr H)=0$: putting 
$$\gamma=\mathsf E(g\xi|\mathscr H) I_{\{|\mathsf E(g\xi|\mathscr H)|\le M\}}\in L^\infty(\mathscr H,D_\xi)$$
and passing in (2.4) to the limit as $M\to\infty$ we conclude that $\mathsf E(g\xi|\mathscr H)=0$ by the monotone convergence theorem. 

The equivalence of the conditions (i) and (iii) for all $p\in [1,\infty]$ follows from the equality $v_p=\|s(a|T_p)\|_q$, which is proved in Lemma 2.4 below. 
\begin{lemma} \label{lem2.3}
Let $\xi\in L^0(\mathscr F,\mathbb R^d)$ and $0\in\ri(\conv \varkappa_\xi)$. If $(\gamma,\xi)\ge 0$ a.s. for some $\gamma\in L^0(\mathscr H,D_\xi)$, then $\gamma=0$ a.s. 
\end{lemma}
{\sc Proof.} Put $A=\{\gamma\neq 0\}$. For any $\omega\in A$ there exists $y\in \varkappa_\xi(\omega)$ such that $(\gamma(\omega),y)<0$ and hence $\int (\gamma(\omega),x)^-\,\mathsf P_\xi(\omega,dx)>0$.  If $\mathsf P(A)>0$ then we obtain the contradiction:
$$ \mathsf E(\gamma ,\xi)^-\ge \mathsf E \mathsf E(I_A (\gamma,\xi)^-|\mathscr H)=\mathsf E\left( I_A \int_{\mathbb R^d} (\gamma(\omega),x)^-\,\mathsf P_\xi(\omega,dx)\right)>0. \ \ \square$$

Lemma 1.1 together with the measurable maximum theorem (Lemma 1.2) imply the existence of an element $h^*_p\in L^0(\mathscr H,T_p)$ such that
$$ s(a(\omega)|T_p(\omega))=(h^*_p(\omega),a(\omega)) \ \textrm{ a.s.}$$

\begin{lemma} \label{lem2.4}
Under the assumptions of Theorem 1.3 we have
$$ v_p=\sup_\gamma\{\mathsf E(\gamma,a): 
 \|(\gamma,\xi)^{-}\|_p\le 1,\ \gamma\in L^\infty(\mathscr H,D_\xi)\}
 =\|s(a|T_p)\|_q, \ \ p\in [1,\infty].$$
\end{lemma}
{\sc Proof.} (a) The case $1\le p<\infty$. Put $U^p_+(\mathscr H)=\{g\in L^p_+(\mathscr H): \|g\|_p \le 1\}$. We have
\begin{eqnarray*}
 U^p_+(\mathscr F) &=& \{g\in L^p_+(\mathscr F): \mathsf E(\mathsf E(g^p|\mathscr H )) \le 1\}\\
&=& \bigcup_{w\in U_+^p(\mathscr H) }\{g\in L^p_+(\mathscr F): \left(\mathsf E(g^p|\mathscr H )\right)^{1/p}\le w\}.
\end{eqnarray*}
Consequently,  
\begin{eqnarray*}  
v_p &=& \sup_\gamma\{\mathsf E(\gamma,a): (\gamma,\xi)^{-}\in U^p_+(\mathscr F),\ \gamma\in L^\infty(\mathscr H,D_\xi)\}\\
&=& \sup_{w\in U_+^p(\mathscr H)} \sup_\gamma\{\mathsf E(\gamma,a): \left(\mathsf E([(\gamma,\xi)^{-}]^p|\mathscr H)\right)^{1/p} \le w,\ \gamma\in L^\infty(\mathscr H,D_\xi)\} 
\end{eqnarray*}
On the set $\{w=0\}$ we have the equality $\mathsf E([(\gamma,\xi)^{-}]^p|\mathscr H)=0$. Therefore,
$$\mathsf E([(\gamma I_{\{w=0\}},\xi)^{-}]^p)=0$$
and $\gamma I_{\{w=0\}}=0$ by Lemma 2.3. Putting $\gamma=w \theta$, where $\theta$ is an $\mathscr H$-measurable vector, we obtain
$$v_p=\sup_{w\in U_+^p(\mathscr H)} \sup_\theta\{\mathsf E w (\theta,a): \mathsf E([(\theta I_{\{w>0\}},\xi)^{-}]^p|\mathscr H)\le 1,\ w \theta\in L^\infty(\mathscr H,D_\xi)\}.$$

Since the values of $\theta$ on the set $\{w=0\}$ do not affect $\mathsf E w (\theta,a)$, by the definition of $T_p$ and the equality $\mathsf E([(\theta,\xi)^{-}]^p|\mathscr H)=\psi_p^p(\omega,\theta(\omega))$ a.s., we get 
$$v_p=\sup_{w\in U_+^p(\mathscr H)} \sup_\theta\{\mathsf E w (\theta,a): \theta\in L^0(\mathscr H,T_p),\ w\theta\in L^\infty(\mathscr H,D_\xi)\}.$$
But $(\theta,a)\le s(a|T_p)$ a.s. for $\theta\in  L^0(\mathscr H, T_p)$. This yields that
\begin{equation} \label{2.5}
v_p\le \sup_{w\in U_+^p(\mathscr H)}\mathsf E(s(a|T_p)w)=
\|s(a|T_p)\|_q.
\end{equation}
We have used Lemma 2.1 in the last equality.

To obtain the inequality, converse to (2.5), put $\theta=h_p^* I_{\{w |h_p^*|\le M\}}$, $M>0$. Clearly, $w\theta\in L^\infty(\mathscr H,D_\xi)$ and
\begin{eqnarray*}
v_p &\ge &\sup_{w\in U_+^p(\mathscr H)}\mathsf E[ w (h_p^*,a) I_{\{w |h_p^*|\le M\}}] =
\sup_{w\in U_+^p(\mathscr H)}\mathsf E(s(a|T_p) w I_{\{w |h_p^*|\le M\}})\\
&=& \|s(a|T_p) I_{\{w |h_p^*|\le M\}} \|_q.
\end{eqnarray*}
By the monotone convergence theorem it follows that $v_p\ge\|s(a|T_p)\|_q$. 

(b) The case $p=\infty$. It follows from
$$ \mathsf P((\gamma,\xi)\ge -1)=\mathsf E \mathsf P (\{\gamma,\xi)\ge -1\}|\mathscr H)=\mathsf E\mathsf P_\xi(\omega,\{x:(\gamma(\omega),x)\ge -1\})$$
that the condition $\|(\gamma,\xi)^{-}\|_\infty\le 1$, meaning that $\mathsf P((\gamma,\xi)\ge -1)=1$, can be represented in the form $\mathsf P_\xi(\omega,\{x:(\gamma,x)\ge -1\})=1$ a.s. In other words, 
$\gamma(\omega)\in -\varkappa_\xi^\circ(\omega)$ a.s.. 

Since $T_\infty=(-\varkappa_\xi^\circ)\cap D_\xi$ this implies that
$$ v_\infty=\sup_\gamma\{\mathsf E (\gamma,a):\gamma\in L^\infty (\mathscr H,(-\varkappa_\xi^\circ)\cap D_\xi)\}\le 
\mathsf E s(a|T_\infty).$$

On the other hand, $h_\infty^* I_{\{|h_\infty^*|\le M\}}\in L^\infty(\mathscr H,(-\varkappa_\xi^\circ)\cap D_\xi)$ for all $M>0$. Therefore,
$$ v_\infty\ge \mathsf E((h_\infty^*,a) I_{\{|h_\infty^*|\le M\}} )=\mathsf E(s(a|T_\infty) I_{\{|h_\infty^*|\le M\}})$$
and $v_\infty\ge \mathsf Es(a|T_\infty)$ by the monotone convergence theorem. $\square$

\section{Proof of Theorem 1 for $p=\infty$}
\setcounter{equation}{0}
As we have already mentioned, Lemma 2.4 yields that conditions (i) and (iii) of Theorem 1.3 are equivalent. Assume that (ii) is satisfied and put $X=(\gamma,\xi)$, $\gamma\in L^\infty(\mathscr H,D_\xi)$. The implication (ii) $\Longrightarrow$ (i) is a consequence of the inequality
\begin{eqnarray}
\mathsf E(f X) &=& \mathsf E(g X)-\mathsf E((g-f)X)\le\mathsf E(\gamma,\mathsf E(g\xi|\mathscr H))+\mathsf E((g-f)X^-)\nonumber\\
&\le & \|g-f\|_1 \|X^-\|_\infty.
\end{eqnarray}

Let us prove that (ii) follows from (iii). We look for $g$ of the form $g=f+\varphi(\omega,\xi(\omega))$, where $\varphi\in L^0_+(\mathscr H\otimes\mathscr B(\mathbb R^d))$. Firstly, the desired function $\varphi$ should satisfy (1.3):
$$ \mathsf E(\varphi\xi|\mathscr H)=\int\varphi(\omega,x)x\,\mathsf P_\xi(\omega,dx)=-a(\omega)\ \ \textrm{a.s.}$$
Secondly, the function $\omega\mapsto \varphi(\omega,\xi(\omega))$ should be $\mathsf P$-integrable. We construct a function $\varphi$ with these properties in Lemma 3.3 after some preliminary work.
\begin{lemma} \label{lem:3.1}
Consider a probability measure $\mathsf Q$ on $(\mathbb R^d,\mathscr B(\mathbb R^d))$ with the support $\varkappa$. If $0\in\ri(\conv\varkappa)$ then for all $y$ in the linear span $D$ of $\varkappa$ the following equality holds true:
$$ w(y):=\inf\left\{\int\varphi(x)\,\mathsf Q(dx):\int\varphi(x) x \,\mathsf Q(dx)=y,\ \varphi\in L^\infty_+(\mathsf Q)\right\}=\mu(y|\conv\varkappa).$$
\end{lemma}
{\sc Proof.} It is easy to check that the epigraph of $w$: $\epi w=\{(y,\alpha)\in D\times\mathbb R: w(y)\le\alpha\}$ is a convex set (see \cite{Rockaf67}, Lemma 2). Following the general scheme of duality theory (see e.g. \cite{Rockaf67}, \cite{MagTik00}) let us find the conjugate function (Young-Fenchel transform) of $w$: 
\begin{eqnarray*}
w^*(\lambda) &=& \sup_{y\in D}\{(y,\lambda)-w(y)\}\\
             &=& \sup_{\varphi,y}\{(y,\lambda)-\int\varphi(x)\mathsf Q(dx):\int\varphi(x)x\,\mathsf Q(dx) = y,\ \varphi\in L^\infty_+(\mathsf Q)\}\\
&=& \sup_{\varphi}\{\int\varphi(x)((x,\lambda)-1)\,\mathsf Q(dx):\varphi\in L^\infty_+(\mathsf Q)\}=\delta(\lambda|\varkappa^\circ),\ \ \lambda\in D.
\end{eqnarray*} 
Here $\delta$ is the indicator function: $\delta(\lambda|\varkappa^\circ)=0$, $\lambda\in\varkappa^\circ$; $\delta(\lambda|\varkappa^\circ)=+\infty$, $\lambda\not\in \varkappa^\circ$. 
The Young-Fenchel transform of $w^*$ is of the form:
$$w^{**}(y) = \sup_{\lambda\in D}\{(y,\lambda)-w^*(\lambda)\}=s(y|\varkappa^\circ)=\mu(y|\conv\varkappa),\ \ y\in D.$$

We claim that $\dom w:=\{y\in D:w(y)<\infty\}=D$. Clearly, this is the case iff the set $A=\{\int\varphi(x) x \,\mathsf Q(dx):\ \varphi\in L^\infty_+(\mathsf Q)\}$ coincides with $D$. 

Assume that $z\in D$ does not belong to the convex set $A$. Then there exists a non-zero vector $h\in D$, separating $A$ and $z$:
$$ \left(\int\varphi(x) x \,\mathsf Q(dx),h\right)=\int\varphi(x) (x,h)\, \mathsf Q(dx)\le (z,h), \ \ \varphi\in L^\infty_+(\mathsf Q).$$
Putting $\varphi(x)=c I_{\{(h,x)\ge 0\}}$, where $c\in\mathbb R_+$, we conclude that the inequality
$$  c \int (x,h)^+\, \mathsf Q(dx)\le (z,h)$$
should hold true for all $c>0$. Consequently $(x,h)^+=0$ $\mathsf Q$-a.s. Then $(h,x)\le 0$, $x\in\varkappa$ and $\varkappa$ is contained in the subspace orthogonal to $h$, since $0\in\ri(\conv\varkappa)$. This means that the linear span of $\varkappa$ does not coincide with $D$, a contradiction. 

Thus, $\dom w=D$, $w$ is continuous on $D$ and $w=w^{**}$ by the Fenchel-Moreau theorem \cite{MagTik00}. $\square$

\begin{lemma} \label{lem:3.2}
There exists a function $\chi:[0,1]\times\mathbb R^d\mapsto\mathbb R$, measurable with respect to $\mathscr B([0,1])\otimes\mathscr B(\mathbb R^d)$ and possessing the following property: for any probability measure $\mathsf Q$ on $\mathscr B(\mathbb R^d)$ and for any $\mathscr B(\mathbb R^d)$-measurable real-valued function $f$ there exists $r\in [0,1]$ such that $\chi(r,x)=f(x)$ $\mathsf Q$-a.s.
\end{lemma}

Lemma 3.2 is borrowed from the paper \cite{EvsShuTak04} (Theorem A.3).

\begin{lemma} \label{lem:3.3} 
If $\xi\in L^1(\mathscr F,\mathbb R^d)$, $0\in\ri(\conv\varkappa_\xi)$ a.s., $a\in L^0(\mathscr H,D_\xi)$ and $\nu=\mu(-a|\conv\varkappa_\xi)$, then there exists a function $\varphi\in L^0_+(\mathscr H\otimes\mathscr B(\mathbb R^d))$ 
such that
$$ \int\varphi(\omega,x)x\,\mathsf P_\xi(\omega,dx)=-a(\omega)\ \textrm{a.s.},$$ 
$$ \int\varphi(\omega,x)\,\mathsf P_\xi(\omega,dx)\in[\nu(\omega),\nu(\omega)+\varepsilon(\omega)]\ \textrm{a.s.} $$
for any $\mathscr H$-measurable function $\varepsilon>0$.
\end{lemma}
{\sc Proof.} Consider the trace $\mathscr H'=\Omega'\cap\mathscr H$ of the $\sigma$-algebra $\mathscr H$ on the set $\Omega'=\{\omega:0\in\ri(\conv\varkappa(\omega))\}\in\mathscr H^\mathsf P$. Let $\chi$ be some function, mentioned in Lemma 3.2. We fix an $\mathscr H$-measurable function $\varepsilon>0$ and introduce the set-valued mapping $G:\Omega'\mapsto [0,1]$ by the formula
\begin{eqnarray*}
G(\omega) &=& \{y\in [0,1]:\int\chi(y,x)\,\mathsf P_\xi(\omega,dx)\in[\nu(\omega),\nu(\omega)+\varepsilon(\omega)],\\ &&\int\chi(y,x)x\,\mathsf P_\xi(\omega,dx)=-a(\omega),
\int\chi^-(y,x)\,\mathsf P_\xi(\omega,dx)=0\}.
\end{eqnarray*}
Applying Lemma 3.1 to $\mathsf Q(dx)=\mathsf P_\xi(\omega,dx)$ and Lemma 3.2, we conclude that  $G(\omega)\neq\varnothing$ for all $\omega\in\Omega'$. The functions 
$$\int\chi^-(y,x)\,\mathsf P_\xi(\omega,dx), \ \ \int\chi(y,x)\,\mathsf P_\xi(\omega,dx),\ \ \int\chi(y,x)x\,\mathsf P_\xi(\omega,dx),$$
depending on $(\omega,y)$, are measurable with respect to $\mathscr H\otimes\mathscr B([0,1])$: see \cite{DalMorWil90}, Lemma 2.2(a). Hence, 
$$\gr G=\{(\omega,y)\in\Omega'\times[0,1]:y\in G(\omega)\}\in\mathscr H'\otimes\mathscr B([0,1])$$ 
and by Aumann's measurable selection theorem there exists an $\mathscr H'$-measurable function $r:\Omega'\mapsto [0,1]$, satisfying the condition $r(\omega)\in G(\omega)$ a.s. on $\Omega'$ (\cite{AliBor07}, Corollary 18.27).
The function $\varphi(\omega,x)=\chi(\widehat r(\omega),x)$, where $\widehat r$ is an $\mathscr H$-measurable modification of $r$, has the desired properties. $\square$

{\sc The end of the proof of Theorem 1.3.} Let us prove that condition (iii) implies (ii) ($p=\infty$). According to the assumption,
 $$s(a|T_\infty)=\mu(-a|\conv\varkappa_\xi)\in L^1(\mathscr H),\ \ \ a=\mathsf E(f\xi|\mathscr H).$$ 
Let $\varepsilon>0$ be some constant. Using the notation of Lemma 3.3, we put $g(\omega)=f(\omega)+\varphi( \omega,\xi(\omega))$. The function $g\ge f$ is $\mathscr F$-measurable, $\mathsf P$-integrable since
\begin{eqnarray}
\mathsf E(\varphi\wedge M) &= & \mathsf E\mathsf E(\varphi\wedge M|\mathscr H)=\mathsf E\int(\varphi(\omega,x)\wedge M)\,\mathsf P_\xi(\omega,dx)\nonumber\\
&\le &\mathsf E\mu(-a|\conv\varkappa_\xi)+\varepsilon,\ \ M>0,
\end{eqnarray}
and satisfies the equality (1.3):
$$ \mathsf E(g\xi|\mathscr H)=a(\omega)+\int\varphi(\omega,x)x\,\mathsf P_\xi(\omega,dx)=0\ \ \textrm{a.s.}
\ \square$$

\section{$N$-period model}
\setcounter{equation}{0}
We turn to $N$-period market model on a filtered probability space, presented in the introductory section. 
In addition to the introduced notation denote by $D_{n-1}(\omega)$ the linear span of $\varkappa_{n-1}(\omega)$.  

Our second main result is the following.
\begin{theorem}
If the process $S_n\in L^\infty(\mathscr F_n,\mathbb R^d)$, $n=0,\dots N$ satisfies the NA property, then for an element $f\in L^1_{++}(\mathscr F,\mathsf P)$ the following conditions are equivalent:
\begin{itemize}
\item[(i)] $v:=\sup\{\mathsf E( fX): \|X^-\|_\infty\le 1,\ X\in K\}<\infty$, where 
$$ K=\{G_N^\gamma:\gamma_n\in L^\infty(\mathscr F_{n-1},D_{n-1}),\ n=1,\dots,N\};$$ 
\item[(ii)] there exist an equivalent to $\mathsf P$ martingale measure $\mathsf Q$, whose density satisfies the inequality $d\mathsf Q/d\mathsf P\ge c f$ with some constant $c>0$;
\item[(iii)] the recurrence relation
\begin{eqnarray*}
\beta_N &=& f,\\
\beta_n &=& \mathsf E(\beta_{n+1}|\mathscr F_n)+\mu(-a_n|\conv\varkappa_n),\ \ 
 a_n=\mathsf E(\beta_{n+1}\Delta S_{n+1}|\mathscr F_n)
\end{eqnarray*} 
specifies the $\mathsf P$-integrable sequence $(\beta_n)_{n=0}^N$.
\end{itemize}
\end{theorem}
{\sc Proof.} (ii) $\Longrightarrow$ (i). This statement follows from an estimate, similar to (3.1). 

(i) $\Longrightarrow$ (iii). Consider the process $X^\gamma=1+G^\gamma$:
$$ X^\gamma_{n+1}=X^\gamma_n+(\gamma_{n+1},\Delta S_{n+1}),\ \ \ X_0^\gamma=1.$$
If the random variable $\beta_n\in L^0_+(\mathscr F_n)$ is well-defined, put
$$ u_n=\sup_\gamma\{\mathsf E (\beta_n X_n^\gamma):X_k^\gamma\ge 0,\ \gamma_k\in L^\infty(\mathscr F_{k-1},D_{k-1}),\ 1\le k\le n\}.$$
By virtue of assumption (i) we have 
$$ u_N\le \sup_\gamma\{\mathsf E (\beta_N X_N^\gamma):X_N^\gamma\ge 0,\ \gamma_k\in L^\infty(\mathscr F_{k-1},D_{k-1}),\ 1\le k\le n\}=\mathsf Ef+v<\infty.$$
 
If $u_{m+1}<\infty$ and the process $\gamma$ satisfies the conditions of the definition of $u_{m+1}$, then $\beta_{m+1}\in L^1(\mathscr F_{m+1})$ and
$$ \mathsf E(\beta_{m+1} X_{m+1}^\gamma)=\mathsf E(X_m^\gamma\mathsf E(\beta_{m+1}|\mathscr F_m))+
\mathsf E(\gamma_{m+1},a_m).$$
Consequently,
\begin{equation}
u_{m+1}\ge \mathsf E(X_m^\gamma\mathsf E(\beta_{m+1}|\mathscr F_m)) + t_{m+1},
\end{equation}
$$ t_{m+1}=\sup_{\gamma_{m+1}}\{\mathsf E(\gamma_{m+1},a_m):X_{m+1}^\gamma\ge 0, \ \gamma_{m+1}\in L^\infty(\mathscr F_m,D_m)\}.$$

The condition $X_{m+1}^\gamma=X_m^\gamma+(\gamma_{m+1},\Delta S_{m+1})\ge 0$ a.s. can be rephrased as 
$$ (\gamma_{m+1}(\omega),x)\ge -X_m^\gamma(\omega),\ \ x\in\varkappa_m(\omega) \ \textrm{a.s.},$$
that is, $\gamma_{m+1}\in - X_m^\gamma\varkappa_m^\circ$ a.s. (see the proof of Lemma 2.4 for $p=\infty$). Here we take into account that $\gamma_{m+1}=0$ a.s., if $(\gamma_{m+1},\Delta S_{m+1})\ge 0$ and $\gamma_{m+1}\in D_m$ a.s. (Lemma 2.3). Thus,
$$ t_{m+1}=\sup_{\gamma_{m+1}}\{\mathsf E(\gamma_{m+1},a_m):\gamma_{m+1}\in L^\infty(\mathscr F_m,-X_m^\gamma\varkappa_m^\circ)\}.$$

The measurability of the set-valued mapping $\varkappa_m^\circ$ with respect to $\mathscr F_m$ follows from
$ \varkappa^\circ_m(\omega)=\bigcap_{i=1}^\infty\{h:(h,\eta_i(\omega))\le 1\},$
where $(\eta_i)_{i=1}^\infty$ is a Castaing representation of $\varkappa_m$ and from Theorem 1M of \cite{Rockaf76}, concerning the measurability of a countable intersection. Owing to the compactness of $\varkappa_m^\circ(\omega)$ a.s., which follows from $0\in\ri(\conv\varkappa_m)$, by the measurable maximum theorem there exists an element $\gamma_{m+1}^*\in L^0(\mathscr F_m,-X_m^\gamma\varkappa_m^\circ)$ such that 
$$(\gamma_{m+1},a_m)\le(\gamma_{m+1}^*,a_m)=s(a_m|-X_m^\gamma \varkappa_m^\circ)=X_m^\gamma \mu(-a_m|\conv\varkappa_m).$$
In particular, $t_{m+1}\le \mathsf E(\gamma_{m+1}^*,a_m)$. On the other hand, by approximation of $\gamma_{m+1}^*$ by the elements $\gamma_{m+1}^* I_{\{|\gamma_{m+1}^*|\le M\}}\in L^\infty(\mathscr F_m,-X_m^\gamma\varkappa^\circ)$, $M\to\infty$, we deduce that
$$\mathsf E(\gamma_{m+1}^*,a_m)=\lim_{M\to\infty}\mathsf E (\gamma_{m+1}^* I_{\{|\gamma_{m+1}^*|\le M\}},a_m)\le t_{m+1}$$
by the monotone convergence theorem. 

By plugging the obtained value $t_{m+1}=\mathsf E\left[X_m^\gamma \mu(-a_m|\conv\varkappa_m)\right]$ in (4.1), we get
$$ v_{m+1}\ge \mathsf E\biggl(\biggl(\mathsf E(\beta_{m+1}|\mathscr F_m)+\mu(-a_m|\conv\varkappa_m)\biggr)X_m^\gamma\biggr)=\mathsf E(\beta_m X_m^\gamma).$$
This inequality holds true under the assumption $X_k^\gamma\ge 0$, $\gamma_k\in L^\infty(\mathscr F_{k-1},D_{k-1})$, $k=1,\dots,m$. Hence, $v_m\le v_{m+1}<\infty$. By induction this implies (iii).

(iii) $\Longrightarrow$ (ii). Put $\nu_n=\mu(-a_n|\conv\varkappa_n)$. Recall that $a_n\in L^0(\mathscr F_n,D_n)$ (see (1.4)). By Lemma 3.3 for any $n=1,\dots,N$ there exists a function $\varphi_n\in L^0_+(\mathscr F_n\otimes\mathscr B(\mathbb R^d))$ such that
\begin{equation}
\int\varphi_n(\omega,x)x\,\mathsf P_n(\omega,dx)=-a_n(\omega)\ \textrm{a.s.},
\end{equation}
\begin{equation}
\int\varphi_n(\omega,x)\,\mathsf P_n(\omega,dx)\in[\nu_n(\omega),\nu_n(\omega)+\beta_n(\omega)]\ \textrm{a.s.}
\end{equation}

Put $\zeta_{n+1}(\omega)=\varphi_n(\omega,\Delta S_{n+1}(\omega))$. The inequality
$$ \mathsf E(\zeta_{n+1}\wedge M)=\mathsf E\int(\varphi_n(\omega,x)\wedge M)\,\mathsf P_n(\omega,dx)\le\mathsf E(\nu_n+\beta_n),$$
similar to (3.2), these functions are $\mathsf P$-integrable. We can rewrite (4.2), (4.3) as follows:
\begin{equation}
\mathsf E(\zeta_{n+1}\Delta S_{n+1}|\mathscr F_n)=-a_n,\ \ \mathsf E(\zeta_{n+1}|\mathscr F_n)=\nu_n+ \varepsilon_n\beta_n,
\end{equation}
where $\varepsilon_n$ is an $\mathscr F_n$-measurable function, taking values in $[0,1]$. 
Put $z_N=1+\zeta_N/f$,
$$z_n=\frac{1}{1+\varepsilon_n}\left(1+\frac{\zeta_n}{\beta_n}\right),\ \ n=1,\dots,N-1;\ \ Z=f\prod_{n=1}^N z_n.$$

We claim that the random variable $Z$ is integrable and
\begin{equation}
\mathsf E(z_{n+1}\dots z_N f|\mathscr F_n)=\beta_n(1+\varepsilon_n),\ \ n=0,\dots N-1.
\end{equation}
By virtue of (4.4) and the definition of $(\beta_n)_{n=0}^N$ we have
\begin{eqnarray*}
\mathsf E(z_N f|\mathscr F_{N-1}) &=& \mathsf E(f|\mathscr F_{N-1})+\mathsf E(\zeta_N|\mathscr F_{N-1})\\
&=& \mathsf E(\beta_N|\mathscr F_{N-1})+\nu_{N-1}+\varepsilon_{N-1}\beta_{N-1}=(1+\varepsilon_{N-1})\beta_{N-1}.
\end{eqnarray*}
Assume that the random variable $z_{m+1}\dots z_N f$ is integrable and (4.5) holds true for $n=m$. 
Then
$$ \mathsf E(I_{\{z_m\le M\}} z_m z_{m+1}\dots z_N f)=\mathsf E(I_{\{z_m\le M\}} z_m\beta_m(1+\varepsilon_m))\le\mathsf E(\beta_m+\zeta_m).$$
Hence, $z_m z_{m+1}\dots z_N f\in L^1(\mathscr F)$. Moreover,
\begin{eqnarray*}
&& \mathsf E(z_m z_{m+1}\dots z_N f|\mathscr F_{m-1}) = \mathsf E(z_m\beta_m(1+\varepsilon_m)|\mathscr F_{m-1})
=\mathsf E(\beta_m+\zeta_m|\mathscr F_{m-1})\\
&=&\mathsf E(\beta_m|\mathscr F_{m-1})+\nu_{m-1}+\varepsilon_{m-1}\beta_{m-1}=(1+\varepsilon_{m-1})\beta_{m-1}.
\end{eqnarray*}
By induction (4.5) hold true for all $n$. In particular, $Z\in L^1(\mathscr F)$.

Consider a probability measure $\mathsf Q$ with the density $d\mathsf Q/d\mathsf P=cZ$, $c=1/\mathsf E Z$. Evidently, $d\mathsf Q/d\mathsf P\ge 2^{-N+1} c f$. Let us check that $\mathsf Q$ is a martingale measure. Put $A_{n-1}\in\mathscr F_{n-1}$. We have
\begin{eqnarray*}
\frac{1}{c}\mathsf E_\mathsf Q(I_{A_{n-1}}\Delta S_n) &=& \mathsf E(\mathsf E(Z|\mathscr F_n) I_{A_{n-1}}\Delta S_n)
=\mathsf E(z_1\dots z_{n}\beta_n(1+\varepsilon_n) I_{A_{n-1}}\Delta S_n)\\
&=&\mathsf E(z_1\dots z_{n-1}I_{A_{n-1}}\mathsf E((\beta_n+\zeta_n)\Delta S_n|\mathscr F_{n-1}))=0
\end{eqnarray*}
since $\mathsf E(\zeta_n\Delta S_n|\mathscr F_{n-1})=-a_{n-1}=-\mathsf E(\beta_n\Delta S_n|\mathscr F_{n-1})$.
$\square$

\section{Examples} 
\setcounter{equation}{0}
In example 5.1 we concretize the formulas of condition (iii) of Theorem 1.3 for a scalar random variable $\xi$ in the case of general probability space. In example 5.2 we consider a one-period model on a countable space. 

Example 5.3 underlines the non-local character of the conditions of Theorem 4.1. Therein we construct a process $(S_0,S_1,S_2)$ with no martingale measure, whose density is bounded from below by a positive constant, but, at the same time, for each of the processes $(S_0,S_1)$, $(S_1,S_2)$ such a measure exists. 

At last, example 5.4 shows that conditions (0.2), (0.3) need not be equivalent for $p=\infty$ even if there exists $z\in L^1_{++}$, satisfying the condition $\mathsf E(X z)=0$, $X\in K$ and the subspace $K$ is generated by a countable collection of elements. 

{\bf Example 5.1.} Consider the case of scalar random variable $\xi$. We use the notation of Theorem 1.3. Assume that $\xi\in L^p(\mathscr F)$, $0\in\ri(\conv\varkappa_\xi)$ and $f\in L^q_+(\mathscr F)$, $1/p+1/q=1$, $p\in [1,\infty]$. 

For $q\in (1,\infty]$ we have
$$ \psi_p(\omega,h)=\int [(hx)^-]^p\,\mathsf P_\xi(\omega,dx)=(h^+)^p\mathsf E((\xi^-)^p|\mathscr H)(\omega)+
(h^-)^p \mathsf E((\xi^+)^p|\mathscr H)(\omega)$$
and condition (iii) shapes to
\begin{eqnarray}
s(a|T_p)&=&\sup_h\{\mathsf E(f\xi|\mathscr H)h:\psi_p(\omega,h)\le 1\}\nonumber\\
&=& \frac{(\mathsf E(f\xi|\mathscr H))^+}{\mathsf E((\xi^-)^p|\mathscr H)^{1/p}}+
\frac{(\mathsf E(f\xi|\mathscr H))^-}{\mathsf E((\xi^+)^p|\mathscr H)^{1/p}}\in L^q(\mathscr H).
\end{eqnarray} 

For $q=1$, $p=\infty$ we have $\conv\varkappa_\xi(\omega)=[\delta_1(\omega),\delta_2(\omega)]$,
$0\in (\delta_1,\delta_2)$ a.s. By virtue of Remark 1.7 condition (iii) becomes
\begin{equation}
\mu(-a|[\delta_1,\delta_2])=\frac{(\mathsf E(f\xi|\mathscr H))^+}{|\delta_1|}+
\frac{(\mathsf E(f\xi|\mathscr H))^-}{\delta_2}\in L^1(\mathscr H).
\end{equation}

{\bf Example 5.2.} Here we slightly generalize the model of \cite{DelSch06} (Remark 6.5.2), \cite{RokSch06} (Example 2). Put $\Omega=\mathbb N$. Consider a countable partition $(A_0^j)_{j=1}^\infty$ of the set $\Omega$:
$$ \mathbb N=\bigcup_{j=1}^\infty A_0^j,\ \ A_0^i\cap A_0^k=\varnothing,\ i\neq k. $$
Denote by $\mathscr H$ the $\sigma$-algebra, generated by this partition.  Let
$$A_0^j=A_1^{2j-1}\cup A_1^{2j},\ \ A_1^{2j-1}\cap A_1^{2j}=\varnothing,\ \ j=1,\dots,\infty$$
and consider the $\sigma$-algebra $\mathscr F$ generated by the sets $(A_1^j)_{j=1}^\infty$. Assume that $\mathsf P(A_1^j)>0$, $j\in\mathbb N$ and let $\xi\in L^p(\mathscr F)$, $1\le p\le\infty$ be a random variable with $0\in\ri(\conv\varkappa_\xi)$:
$$\xi(\omega)>0,\ \omega\in A_1^{2j-1},\ \ \ \xi(\omega)<0,\ \omega\in A_1^{2j},\ \ j\in\mathbb N.$$

Let $f\in L^q_+(\mathscr F)$, $1/p+1/q=1$, $p\in [0,\infty]$. For brevity, we put $\eta^j=\eta(\omega)$, $\omega\in A_1^j$ for any $\mathscr F$-measurable random variable $\eta$. Define the random variable $\rho$ by the formula 
$$ \rho(\omega)=\sum_{j=1}^\infty\left(f^{2j}\left|\frac{\xi^{2j}}{\xi^{2j-1}}\right|\frac{\mathsf P(A_1^{2j})}{\mathsf P(A_1^{2j-1})} I_{A_1^{2j-1}}(\omega)+ f^{2j-1}\left|\frac{\xi^{2j-1}}{\xi^{2j}}\right|\frac{\mathsf P(A_1^{2j-1})}{\mathsf P(A_1^{2j})} I_{A_1^{2j}}(\omega)\right).$$

We claim that a necessary and sufficient condition for the existence of a random variable $g$, satisfying conditions (ii) of Theorem 1.3, is the following: 
\begin{equation}
\rho\in L^q(\mathscr F).
\end{equation}

We make use of conditions (5.1), (5.2), obtained in example 5.1. In our case
$$ \mathsf E(f\xi|\mathscr H)(\omega)=\sum_{j=1}^\infty \frac{f^{2j-1}\xi^{2j-1}\mathsf P(A_1^{2j-1})+
f^{2j}\xi^{2j}\mathsf P(A_1^{2j})}{\mathsf P(A_0^j)}I_{A_0^j}(\omega),$$
$$(\mathsf E(f\xi|\mathscr H))^+(\omega)=\sum_{j=1}^\infty \frac{|\xi^{2j}|\mathsf P(A_1^{2j})(\rho^{2j}-f^{2j})^+}{\mathsf P(A_0^j)}I_{A_0^j}(\omega),$$
$$(\mathsf E(f\xi|\mathscr H))^-(\omega)=\sum_{j=1}^\infty \frac{\xi^{2j-1}\mathsf P(A_1^{2j-1})(f^{2j-1}-\rho^{2j-1})^-}{\mathsf P(A_0^j)}I_{A_0^j}(\omega).$$

Let $q=1$. Since $[\delta_1,\delta_2]=\sum_{j=1}^\infty [\xi^{2j},\xi^{2j-1}]I_{A_0^j}$, condition (5.2) shapes to
\begin{eqnarray*}
\mathsf E \mu(-a|[\delta_1,\delta_2]) &=& \sum_{j=1}^\infty\biggl((\rho^{2j}-f^{2j})^+\mathsf P(A_1^{2j})
+(f^{2j-1)}-\rho^{2j-1})^- \mathsf P(A_1^{2j-1})\biggr)\\ 
&=&\|(\rho-f)^+\|_1<\infty,
\end{eqnarray*}
which is equivalent to (5.3), as long as $f\in L^1(\mathscr F)$.

For $q\in (1,\infty]$ we use (5.1). By virtue of the equalities
$$ \mathsf E((\xi^-)^p|\mathscr H)=\sum_{j=1}^\infty\frac{|\xi^{2j}|^p \mathsf P(A_1^{2j})}{\mathsf P(A_0^j)} I_{A_0^j},\ \ \
  \mathsf E((\xi^+)^p|\mathscr H)=\sum_{j=1}^\infty\frac{(\xi^{2j-1})^p \mathsf P(A_1^{2j-1})}{\mathsf P(A_0^j)} I_{A_0^j}$$ 
we get 
$$s(a|T_p)=\sum_{j=1}^\infty\frac{(\rho^{2j}-f^{2j})^+\mathsf P(A_1^{2j})^{1-1/p}
+(f^{2j-1}-\rho^{2j-1})^- \mathsf P(A_1^{2j-1})^{1-1/p}}{\mathsf P(A_0^j)^{1-1/p}} I_{A_0^j}.$$
For $q\in (1,\infty)$ condition (5.1) means that
\begin{eqnarray*}
\|s(a|T_p)\|_q^q &=& \sum_{j=1}^\infty\biggl( [(\rho^{2j}-f^{2j})^+]^q \mathsf P(A_1^{2j})+
+[(f^{2j-1}-\rho^{2j-1})^-]^q \mathsf P(A_1^{2j-1})\biggr)\\
&=&\|(\rho-f)^+\|_q^q<\infty,
\end{eqnarray*}
and is reduced to (5.3). At last, condition $s(a|T_1)\in L^\infty(\mathscr H)$ for $f\in L^\infty(\mathscr F)$ is equivalent to the boundedness of $\rho$. 

{\bf Example 5.3.} Put $\Omega=\mathbb N$ and consider the filtration $\mathscr F_0\subset\mathscr F_1\subset\mathscr F_2$, where the $\sigma$-algebra $\mathscr F_n$ is generated by the sets $(A_n^j)_{j=1}^\infty$, $n=0,1,2$,
$$ A_0^j=\{4j-3,4j-2,4j-1,4j\},\ \ A_1^j=\{2j-1,2j\},\ \ A_2^j=\{j\}.$$
Define the probability measure $\mathsf P$ on $\mathscr F_2=\mathscr F$ by $\mathsf P(A_2^{2j-1})=\mathsf P(A_2^{2j})=2^{-j-1}$.
Note that
$$ \mathsf P(A_1^j)=\mathsf P(A_2^{2j-1})+\mathsf P(A_2^{2j})=2^{-j},$$
$$ \mathsf P(A_0^j)=\mathsf P(A_1^{2j-1})+\mathsf P(A_1^{2j})=2^{-2j+1}+2^{-2j}=\frac{3}{2^{2j}}.$$
We put
$$\xi_1(\omega)=\Delta S_1(\omega)=\sum_{j=1}^\infty\left(I_{A_1^{2j-1}}(\omega)-\frac{1}{2^j}I_{A_1^{2j}}(\omega)\right),$$
$$\xi_2(\omega)=\Delta S_2(\omega)=\sum_{j=1}^\infty\left(I_{A_2^{2j-1}}(\omega)-\frac{1}{2^{j/2}}I_{A_2^{2j}}(\omega)\right).$$

According to Example 5.2, for the existence of $g_n\in L^1(\mathscr F_n)$, $n=1,2$, satisfying the conditions
$$ \mathsf E(g_n\xi_n|\mathscr F_{n-1})=0,\ \ \ g_n\ge 1,$$
it is necessary and sufficient that the functions 
$$\rho_n(\omega)=\sum_{j=1}^\infty\left(\frac{1}{2^{j/n}} \frac{\mathsf P(A_n^{2j})}{\mathsf P(A_n^{2j-1})} I_{A_n^{2j-1}}(\omega)
+ 2^{j/n}\frac{\mathsf P(A_n^{2j-1})}{\mathsf P(A_n^{2j})}I_{A_n^{2j}}(\omega)\right),\ \ n=1,2$$
in the conditions of the form (5.3), are integrable.
A simple calculation shows that it is the case: 
$$\mathsf E \rho_1=\mathsf E\sum_{j=1}^\infty\left(\frac{1}{2^{j+1}} I_{A_1^{2j-1}}+2^{j+1} I_{A_1^{2j}}\right)
=\sum_{j=1}^\infty\left(\frac{1}{2^{3j}}+\frac{2}{2^j}\right)<\infty,$$
$$\mathsf E \rho_2=\mathsf E\sum_{j=1}^\infty\left(\frac{1}{2^{j/2}}I_{A_2^{2j-1}}+2^{j/2} I_{A_2^{2j}}\right)
  =\sum_{j=1}^\infty\left(\frac{1}{2^{3j/2+1}}+\frac{1}{2^{j/2+1}}\right)<\infty.$$

Nevertheless, as we shall see, in the two-period model under consideration, there is no equivalent martingale measure $\mathsf Q$ with the density $d\mathsf Q/d\mathsf P\ge c>0$, where $c$ is some constant. 

Let $\omega\in A_1^j$. With the notation of Theorem 4.1 we have $\beta_2=1$,
$$ a_1(\omega)=\mathsf E(\xi_2|\mathscr F_1)(\omega)=\frac{\mathsf E(\xi_2 I_{A_1^j})}{\mathsf P(A_1^j)}=
\frac{\mathsf P(A_2^{2j-1})-2^{-j/2} \mathsf P(A_2^{2j})}{\mathsf P(A_1^j)}=\frac{1-2^{-j/2}}{2},$$
$$\mu(-a_1|\conv\varkappa_1)(\omega)=\inf\{\lambda>0:-a_1(\omega)\in\lambda [-2^{-j/2},1]\}=2^{j/2} a_1(\omega),$$
$$\beta_1(\omega)=1+\mu(-a_1|\conv\varkappa_1)(\omega)=1+2^{j/2}\frac{(1-2^{-j/2})}{2}=\frac{2^{j/2}+1}{2}$$
and $\mathsf E\beta_1=\sum_{j=1}^\infty (2^{j/2}+1)\mathsf P(A_1^j)/2<\infty$.

Now assume that $\omega\in A_0^j$. Then
\begin{eqnarray*}
a_0(\omega)\mathsf P(A_0^j)&=&\mathsf E(\beta_1\xi_1|\mathscr F_0)(\omega)\mathsf P(A_0^j)=\mathsf E(\beta_1\xi_1 I_{A_0^j})=\frac{2^{j-1/2}+1}{2} \mathsf P(A_1^{2j-1})\\
 &-& \frac{2^j+1}{2}\frac{1}{2^j} \mathsf P(A_1^{2j})=
\frac{1}{2^{2j}}\left(2^{j-1/2}+\frac{1}{2}-\frac{1}{2^{j+1}}\right)
\end{eqnarray*}
In addition, $a_0(\omega)>0$ and 
$$\mu(-a_0|\conv\varkappa_0)(\omega)=\inf\{\lambda>0:-a_0(\omega)\in\lambda [-2^{-j},1]\}=2^j a_0(\omega).$$
This yields that 
$$\mathsf E \mu(-a_0|\conv\varkappa_0)=\sum_{j=1}^\infty 2^j a_0^j\mathsf P(A_0^j)=\infty,\ \ a_0^j=a_0(\omega),\ \omega\in A_0^j.$$
Therefore,
$\beta_0=\mathsf E(\beta_1|\mathscr F_0)+\mu(-a_0|\conv\varkappa_0)\not\in L^1(\mathscr F_0).$

This result shows also that
$$ \sup_\gamma\{\mathsf E G_2^\gamma: \gamma_n\in L^\infty(\mathscr F_{n-1}),\ n=1,2,\ G_2^\gamma\ge -1\}=\infty,$$
whereas
$$ \sup_{\gamma_n}\{\mathsf E (\gamma_n,\xi_n): \gamma_n\in L^\infty(\mathscr F_{n-1}),\ (\gamma_n,\xi_n)\ge -1\}<\infty,\ \ n=1,2.$$

Let us present a strategy $\gamma_n\in L^0(\mathscr F_{n-1})$, $n=1,2$, satisfying the conditions 
$$ \mathsf E G_2^\gamma=\infty,\ \ \ G_2^\gamma\ge -1.$$ 
The strategy, constructed below, is ''aggressive'' and consists in buying of the maximal allowable amount of stocks in each step. 

Put $\gamma_1(\omega)=\sum_{j=1}^\infty 2^j I_{A_0^j}$. Then
$$G_1^\gamma=\sum_{j=1}^\infty\left(2^j I_{A_1^{2j-1}}-I_{A_1^{2j}}\right)\ge -1.$$
Since $A_1^{2j-1}=A_2^{4j-3}\cup A_2^{4j-2}$ and 
$$\xi_2(\omega)=1,\ \ \omega\in A_2^{4j-3},\ \ \ \xi_2(\omega)=-\frac{1}{2^{j-1/2}},\ \ \omega\in A_2^{4j-2},$$
we see that the portfolio
$\gamma_2(\omega)=\sum_{j=1}^\infty 2^{j-1/2} (2^j+1) I_{A_1^{2j-1}}$
is admissible:
$$G_2^\gamma=\sum_{j=1}^\infty\left(2^j I_{A_1^{2j-1}}-I_{A_1^{2j}}\right)+
\sum_{j=1}^\infty \left( 2^{j-1/2} (2^j+1) I_{A_2^{4j-3}}-(2^j+1) I_{A_2^{4j-2}}\right)\ge -1$$
and $\mathsf E G_2^\gamma=\infty$ as long as
$$ \mathsf P(A_1^{2j-1})=2^{-2j+1},\ \ \mathsf P(A_1^{2j})=2^{-2j},\ \
   \mathsf P(A_2^{4j-3})=\mathsf P(A_2^{4j-2})=2^{-2j}.$$
   
{\bf Example 5.4.} Let $\Omega=\mathbb N$, $\mathscr F_0=\{\varnothing,\Omega\}$ and let $\mathscr F$ be generated by one-point subsets of $\mathbb N$. We put $A_j=\bigcup_{i=j}^\infty\{2i\}$, $B_j=\{4j+1\}$, 
$$ \Delta S_1^j=\xi^j=2^j I_{B_{j-1}}-I_{A_j}, \ \ j\in\mathbb N$$  
and define the probability measure  $\mathsf Q$ on $\mathscr F$ by $\mathsf Q\{2j-1\}=\mathsf Q\{2j\}=2^{-j-1}$. Clearly, $\mathsf Q$ is a martingale measure for $S$: 
$$\mathsf Q(B_{j-1})=\mathsf Q\{2(2j-1)-1\}=\frac{1}{2^{2j}},\ \ \ 
  \mathsf Q(A_j)=\sum_{i=j}^\infty\frac{1}{2^{i+1}}=\frac{1}{2^j}$$
$$\mathsf E_\mathsf Q\xi^j=2^j\mathsf Q(B_{j-1})-\mathsf Q(A_j)=0.$$

Put $B=\cup_{j=1}^\infty B_{j-1}$ and $B'=\Omega\backslash(A_1\cup B)=\cup_{j=1}^\infty\{4j-1\}$. The set $\Omega$ coincides with the union of disjoint sets $A_1$, $B$, $B'$. We note that
$$ \mathsf Q(A_1)=\frac{1}{2}, \ \ \ \mathsf Q(B')=\sum_{j=1}^\infty\mathsf Q\{2(2j)-1\}=\sum_{j=1}^\infty\frac{1}{2^{2j+1}}=\frac{1}{6}$$
and define an equivalent to $\mathsf Q$ ''market'' measure $\mathsf P$ by 
$$ \mathsf P(C)=\mathsf E_\mathsf Q(\zeta I_C),\ \ \ \zeta=\sum_{i=1}^\infty 2^{i-1} I_{B_{i-1}}+\frac{3}{4}(I_{A_1}+I_{B'}),\ \ C\in\mathscr F.$$  

Let $J$ be a finite subset of $N$. Putting in the inequality 
\begin{equation}
G^\gamma(\omega):=\sum_{j\in J}\gamma^j\xi^j(\omega)\ge -1,\ \ \ \omega\in\mathbb N
\end{equation}
$\omega=2m>\max J$ and then $\omega=4(m-1)+1$, we get:
$$ \sum_{j\in J}\gamma^j\le 1,\ \ \ 2^m\gamma^m\ge -1.$$  
As far as
$$ \mathsf E_\mathsf Q (\zeta\xi^j)=\mathsf E_\mathsf Q (2^{2j-1} I_{B_{j-1}}-\frac{3}{4} I_{A_j})=\frac{1}{2}-\frac{3}{4}\frac{1}{2^j},$$
for $\gamma$ satisfying (5.4) we have
$$ \mathsf E G^\gamma=\sum_{j\in J}\gamma^j \mathsf E_\mathsf Q (\zeta\xi^j)
=\frac{1}{2}\sum_{j\in J}\gamma^j -\frac{3}{4}\sum_{j\in J}\gamma^j 2^{-j}\le \frac{1}{2}+\frac{3}{4}\sum_{j=1}^\infty\frac{1}{2^{2j}}=\frac{3}{4}.$$
   
On the other hand, if $g$ is the $\mathsf P$-density of a martingale measure and $g$ is uniformly bounded from below by a constant $c>0$, then
$$ \mathsf E(g\xi^j)=2^j\mathsf E (g I_{B_{j-1}})-\mathsf E(g I_{A_j})=0,$$   
$$ \mathsf E(g I_{A_j})\ge c 2^j \mathsf P(B_{j-1})=c 2^{2j-1}\mathsf Q (B_{j-1})=\frac{c}{2},$$
in contradiction to the dominated convergence theorem, since $\lim_{j\to\infty} I_{A_j}=0$ a.s. 

Summing up, for the subspace $K\subset L^\infty(\mathscr F)$, generated by the countable collection of elements $(\xi^j)_{j=1}^\infty$, condition (0.3) is satisfied for $f=1$, $p=\infty$. Moreover, there exists and element $z=\zeta^{-1}\in L^1_{++}(\mathscr F)$ such that $\langle X,z\rangle=\mathsf E(X z)=\mathsf E_\mathsf Q X=0$, $X\in K$. However, there is no element $g$, satisfying (0.2) for $q=1$:
a counterexample to the assertion of Theorem 6.1 of \cite{Lei08}.


\begin{thebibliography}{11}
\bibitem{AliBor07}
{\it Aliprantis C.D., Border K.C.} Infinite dimensional analysis. A hitchhicker's guide. 3rd edn. Berlin:\,Springer, 2007, 703 p.
\bibitem{Bur68}
{\it Bourbaki N.} General Topology. Fundamental Structures, Moscow:\,Nauka, 1968, 272 p, [Russian translation].
\bibitem{DalMorWil90}
{\it Dalang R.C., Morton A., Willinger W.} Equivalent martingale measures and no-arbitrage in stochastic securities market models. --- Stoch. Stoch. Rep., 1990, v. 29, No 2, p. 185--201.
\bibitem{DelSch06}
{\it Delbaen F., Schachermayer W.} The mathematics of arbitrage. Berlin:\,Springer, 2006, 373 p.
\bibitem{EvsShuTak04}
{\it Evstigneev I.V., Sch\"{u}rger K., Taksar M.I.} On the fundamental theorem of asset pricing: random constraints and bang-bang no-arbitrage criteria. --- Math. Finance, 2004, v. 14, No. 2, p. 201--221.
\bibitem{Lei08}
{\it Leitner J.} Optimal portfolios with lower partial moment constraints and LPM-risk-optimal martingale measures. --- Math. Finance, 2008, v. 18, No. 2, p. 317–331.
\bibitem{MagTik00}
{\it Magaril-Il'yaev G.G., Tikhomirov V.M.} Convex analyis: theory and applications. Providence, RI: AMS, 2003, 183 p.
\bibitem{RasSte05}
{\it R\'asonyi M., Stettner L.} On utility maximization in discrete-time financial market models. 
---  Ann. Appl. Probab., 2005, v. 15, No. 2, p. 1367--1395.
\bibitem{Rockaf67}
{\it Rockafellar R.T.} Duality and stability in extremum problems involving convex functions. --- Pacific J. Math., 1967, v. 21, No. 1, p. 167--187. 
\bibitem{Rockaf76}
{\it Rockafellar R.T.} Integral functionals, normal integrands and measurable selections. --- Lecture Notes in Math., 1976, v. 543, p. 157--207.
\bibitem{Rok05}
{\it Rokhlin  D.B.} The Kreps-Yan theorem for $L^\infty$. --- Int. J. Math. Math. Sci., 2005, v. 2005, No. 17, p. 2749--2756.
\bibitem{RokSch06}
{\it Rokhlin  D., Schachermayer W.} A note on lower bounds
of martingale measure densities. --- Illinois J. Math., 2006, v. 50, No. 4, p. 815--824.
\bibitem{Shi98}
{\it Shiryaev A.N.} Essentials of Stochastic Finance. Vol.2: Theory, Moscow:\,Phasis, 1998, 528 p., [Russian edition]
\end{thebibliography}
\end{document}